\documentclass[10pt]{amsart}

\usepackage{geometry}
\usepackage[all]{xy}
\usepackage{graphicx}
\usepackage{amssymb}
\usepackage{amsfonts}
\usepackage{epstopdf}
\usepackage{amsmath, amsthm}
\usepackage{mathrsfs}
\usepackage{rotating}

\newcommand{\bb}[1]{\mathbb{#1}}
\newcommand{\cc}[1]{\mathcal{#1}}
\newcommand{\lie}[1]{\mathfrak{#1}}

\newcommand{\diag}{\textrm{diag}}
\def\inv{^{-1}}
\newcommand{\kt}{\Box (\pi) \times \triangle (r)}

\theoremstyle{plain}
\newtheorem{theorem}{Theorem}[section]

\newtheorem{corollary}[theorem]{Corollary}
\newtheorem{lemma}[theorem]{Lemma}
\newtheorem{proposition}[theorem]{Proposition}

\theoremstyle{definition}
\newtheorem{example}[theorem]{Example}
\newtheorem{remark}[theorem]{Remark}

\begin{document}
\title{Gromov width of non-regular coadjoint orbits of $U(n)$, $SO(2n)$ and $SO(2n+1)$.}
\author{Milena Pabiniak}

\address{Milena Pabiniak,
CAMGSD, Departamento de Matem\'atica,
Instituto Superior T\'ecnico, Lisboa, Portugal}
\email{mpabiniak@math.ist.utl.pt}
\thanks{The author was supported by the Funda\c{c}\~ao para a Ci\^encia e a Tecnologia (FCT, Portugal) grant SFRH/BPD/87791/2012 during the final stage of this research.}

\maketitle 
\begin{abstract}Let $G$ be a compact connected Lie group $G$ and $T$ its maximal torus. The coadjoint orbit $\cc{O}_{\lambda}$ through $\lambda \in \lie{t}^*$ is canonically a symplectic manifold. Therefore we can ask the question about its Gromov width.
In many known cases the Gromov width is exactly the minimum 
over the set $\{ \langle \alpha_j^{\vee},\lambda \rangle; \alpha_j^{\vee} \textrm{ a coroot, }\langle \alpha_j^{\vee},\lambda \rangle>0\}$. 
We show that the Gromov width of coadjoint orbits of the unitary group and of most of the coadjoint orbits of the special orthogonal group is at least the above minimum. The proof uses the torus action coming from the Gelfand-Tsetlin system.\end{abstract}
\tableofcontents

\section{Introduction}
Let $(M,\omega)$ be a symplectic manifold. Non-degeneracy of $\omega$ implies that every symplectomorphism is a volume preserving transformation. However,
Gromov's Non-squeezing theorem proves that a group of symplectomorphisms is a proper subset of the group of volume preserving transformations. The theorem says that a ball
$B^{2N}(r)$ of radius $r$, in a symplectic vector space $\bb{R}^{2N}$ with the usual symplectic structure, can be symplectically embedded into $B^2(R)\times \bb{R}^{2N-2}$ only if $r \leq R$.
This motivated the definition of the invariant called {\it Gromov width.}
Consider a ball of {\bf capacity} $a$
$$ B^{2N}_a = \Big \{ z \in \bb{C}^N \ \Big | \ \pi \sum_{i=1}^N |z_i|^2 < a \Big \} , $$
with the standard symplectic form
$\omega_{std} = \sum dx_j \wedge dy_j$.
The \textbf{Gromov width} of a $2N$-dimensional symplectic manifold $(M,\omega)$
is the supremum of the set of $a$'s such that $B^{2N}_a$ can be symplectically
embedded in $(M,\omega)$. 

In this work we focus on the Gromov width of  coadjoint orbits of Lie groups. 
A Lie group $G$ acts on itself by conjugation
$$G \ni g:G \rightarrow G,\;\;\;g(h)=g h g \inv.$$
Derivative of the above map taken at the identity element gives the action of $G$ on its Lie algebra $\lie{g}$, called the adjoint action.
This induces the action of $G$ on $\lie{g}^*$, the dual of its Lie algebra, called the coadjoint action. Each orbit 
$\mathcal{O}$ of the coadjoint action 
is naturally equipped with the Kostant-Kirillov symplectic form:
$$\omega_{\xi}(X,Y)=\langle \xi, [X,Y]\rangle,\;\;\;\xi \in \lie{g}^*,\;X,Y \in \lie{g}.$$
For example, when $G=U(n)$
the group of (complex) unitary matrices, a coadjoint orbit can be identified with the set of Hermitian matrices with
a fixed set of eigenvalues. With this identification, the coadjoint action of $G$ on an orbit $\mathcal{O}$ is simply the action by conjugation. It is Hamiltonian, and the momentum map is just inclusion $\mathcal{O}\hookrightarrow\lie{g}^*$. We recall the notions of Hamiltonian actions and momentum maps in Section \ref{tech}.

Choose a maximal torus $T \subset G$ and a positive Weyl chamber $\lie{t}^*_+$. Every coadjoint orbit intersects the positive Weyl chamber in a single point. Therefore there is a bijection between the coadjoint orbits and points in the positive Weyl chamber. Points in the interior of the positive Weyl chamber are called {\bf regular} points.
The main result of this paper describes a lower bound for Gromov width of the coadjoint orbits of the unitary group.
\begin{theorem}\label{width}
Let $M:=\mathcal{O}_{\lambda}$ be the coadjoint orbit of $G=U(n)$ through a point $\lambda \in (\lie{t}_G)^*_+$ (regular or not) or of $G=SO(2n+1),\, SO(2n)$ through a point $\lambda \in (\lie{t}_G)^*_+$ satisfying condition $(*)$ stated below.
The Gromov width of $M$ is at least
$$r_{G}(\lambda):=\min\{\, \left\langle \alpha^{\vee},\lambda \right\rangle \,; \alpha^{\vee} \textrm{ a  coroot and } \left\langle \alpha^{\vee},\lambda \right\rangle>0\}.$$
\end{theorem}
To state the condition $(*)$ we need to review the root system of the special orthogonal groups and fix the notation. Therefore we delay the explanation of $(*)$ till Theorem \ref{mainso}. Here we only note that all regular orbits satisfy condition $(*)$.

This particular lower bound is important because in many known cases it describes the Gromov width, not only its lower bound. 
Karshon and Tolman in \cite{KT1} showed that the Gromov width of complex Grassmannians is given by the above formula. Zoghi in \cite{Z} analyzed orbits satisfying some additional integrality conditions. He called an orbit $\cc{O}_{\lambda}$ {\bf indecomposable} if there exists a simple root $\alpha$ such that for each root $\alpha'$ there exists a positive integer $k$ (depending on $\alpha'$) such that
$$k\,\langle \alpha^{\vee},\lambda \rangle=\langle (\alpha')^{\vee},\lambda \rangle .$$ 
\begin{theorem}\cite[Proposition 3.16]{Z}
For compact connected simple Lie group $G$ the formula 
$\min\{\, \left|\left\langle \alpha^{\vee},\lambda \right\rangle \right|\,; \alpha^{\vee} \textrm{ a  coroot}\}$
gives an upper bound for Gromov width of regular indecomposable $G$-coadjoint orbit through $\lambda$.
\end{theorem}
Combinining these results we obtain
\begin{theorem}
The Gromov width of a regular indecomposable $U(n)$- or $SO(n)$-coadjoint orbit through $\lambda$ is exactly $$\min\{\, \left|\left\langle \alpha^{\vee},\lambda \right\rangle \right|\,; \alpha^{\vee} \textrm{ a  coroot}\}.$$
\end{theorem}
The result of Zoghi was recently extended by Caviedes in \cite{Cav} to some non-regular $U(n)$ orbits. We quote his result in the $U(n)$ subsection below.
\subsection{Reformulation of the main result for the unitary group.}
Choose as the maximal torus $T$ of $U(n)$ a the subgroup of diagonal matrices. We use the following indentifications: 
\\ - the exponential map $\exp \colon Lie(S^1) \rightarrow S^1$ is given by $t \rightarrow e^{2 \pi i t}$,
\\ - $\lie{u}(n)$ is identified with the set of $n \times n$ Hermitian matrices, 
\\ - the pairing in $\lie{u}(n)$,
$(A,B)=\textrm{trace}(AB)$
gives us the identification of $\lie{u}^*(n)$ with $\lie{u}(n)$,
\\ - $\lie{t}^*$ and  $\lie{t}$ are identified with diagonal Hermitian matrices, and then with $\bb{R}^n$ (by mapping a diagonal matrix to its diagonal entries);
\\Kernel of the exponential maps forms a lattice in  $\lie{t}$ and thus induces a lattice in $\lie{t}^*$. Choose the following chamber
$$(\lie{t}^*)_+:=\{ ( \lambda_{11}, \lambda_{22}, \ldots , \lambda_{nn});\,  \lambda_{11} \geq \lambda_{22} \geq \ldots  \geq \lambda_{nn}\}$$
to be the positive Weyl chamber. Fix any $\lambda \in (\lie{t}^*)_+$, regular or not and denote by $\cc{O}_{\lambda}$ the $U(n)$-coadjoint orbit through $\lambda$.
Recall that the root system of $U(n)$ consists of vectors $\pm (e_j-e_k)$, $j\neq k$, of lattice length $2$. The pairing of $\lambda$ with a coroot $(e_{j}-e_{k})^{\vee}$ gives 
$$\, \left\langle   (e_{j}-e_{k})^{\vee} , \lambda \right\rangle \,= 2\, \frac{\left\langle e_{j}-e_{k},\lambda\right\rangle}{\left\langle e_{j}-e_{k},e_{j}-e_{k}\right\rangle}= (\lambda_j-\lambda_k).$$ 
Note that the real dimension of a $U(n)$ coadjoint orbit $\cc{O}_{\lambda}$ through a point $\lambda$ with $\lambda_1=\ldots = \lambda_{n_1} >\lambda_{n_1+1}=\ldots =\lambda_{n_2}> \ldots >\lambda_{n_{m}+1}=\ldots =\lambda_n$
is $$\dim_{\bb{R}}(\cc{O}_{\lambda})=n(n-1)-(\,k_1(k_1-1) + \ldots + k_{m+1}(k_{m+1}-1) \,)$$
where $k_j=n_{j+1}-n_j$, $n_{m+1}=n$.

Now we restate the main theorem for the unitary group in more explicit form.

\begin{theorem}\label{mainun}
 Let $\lambda=(\lambda_1, \ldots, \lambda_n) \in (\lie{t}^*)_+$  and let $m,n_1,\ldots,n_m$ be integers such that 
$$\lambda_1=\ldots = \lambda_{n_1} >\lambda_{n_1+1}=\ldots =\lambda_{n_2}> \ldots >\lambda_{n_{m}+1}=\ldots =\lambda_n.$$
 The Gromov width of $\cc{O}_{\lambda}$, $U(n)$ coadjoint orbit through $\lambda$, is at least 
$$r_{U(n)}(\lambda):=\min\{\,\lambda_{n_1}-\lambda_{n_1+1}, \lambda_{n_2}-\lambda_{n_2+1},\ldots, \lambda_{n_m}-\lambda_{n_m+1} \}.$$
\end{theorem}
Caviedes in \cite{Cav} proved the following result.
\begin{theorem}\cite[Theorem 5.4]{Cav}
Let $\lambda=(\lambda_1, \ldots, \lambda_n) \in (\lie{t}^*)_+$ and suppose that there are indicies $i,j\in \{1,\ldots,n\}$ such that for any $i',j'\in \{1,\ldots,n\}$ the difference $\lambda_{i'}-\lambda_{j'}$ is an integer multiple of $\lambda_i-\lambda_j$. Then the Gromov width of $\cc{O}_{\lambda}$ is at most $|\lambda_i-\lambda_j|.$
\end{theorem}
Note that in that case we have $r_{U(n)}(\lambda)=|\lambda_i-\lambda_j|.$ Combining these results together we can calculate the actual Gromov width.
\begin{theorem}
Let $\lambda=(\lambda_1, \ldots, \lambda_n) \in (\lie{t}^*)_+$ and suppose that there are indicies $i,j\in \{1,\ldots,n\}$ such that for any $i',j'\in \{1,\ldots,n\}$ the difference $\lambda_{i'}-\lambda_{j'}$ is an integer multiple of $\lambda_i-\lambda_j$. Then the Gromov width of $\cc{O}_{\lambda}$ is exactly $$|\lambda_i-\lambda_j|.$$
\end{theorem}
\subsection{Reformulation of the main result for the special orthogonal group.}
We identify the Lie algebra $\lie{so(m)}$, and its dual $\lie{so(m)}^*$ with the vector space of skew symmetric matrices of appropriate size.
Throughout the paper we use the notation
\begin{displaymath}R(\alpha)=
\left(\begin{array}{cc}
    \cos(\alpha)&  -\sin(\alpha)\\
\sin(\alpha)& \cos(\alpha)
    \end{array}
\right),\,\,\,
L(a)=
\left(\begin{array}{cc}
    0&  -a\\
a& 0
    \end{array}
\right)
\end{displaymath} 
and make the following choices of maximal tori
\scriptsize\begin{displaymath} T_{SO(2n+1)}=\left\{ \left( \begin{array}{ccccc}
R(\alpha_1) &&&&\\
&R(\alpha_2)&&&\\
&& \ddots &&\\
&&& R(\alpha_n)&\\
&&&&1
\end{array}\right) \right\},\,\,\,T_{SO(2n)}=\left\{\left( \begin{array}{cccc}
R(\alpha_1) &&&\\
&R(\alpha_2)&&\\
&& \ddots &\\
&&& R(\alpha_n)
\end{array}\right)\right\} \end{displaymath}
\normalsize where $\alpha_j\in S^1.$  The corresponding Lie algebra duals are \scriptsize\begin{displaymath} \lie{t}_{SO(2n+1)}^*=\left\{ \left( \begin{array}{ccccc}
L(a_1) &&&&\\
&L(a_2)&&&\\
&& \ddots &&\\
&&& L(a_n)&\\
&&&&0
\end{array}\right) \right\},\,\,\,\lie{t}_{SO(2n)}^*=\left\{\left( \begin{array}{cccc}
L(a_1) &&&\\
&L(a_2)&&\\
&& \ddots &\\
&&& L(a_n)
\end{array}\right)\right\}. \end{displaymath}
\normalsize 
We identify these duals with $\bb{R}^n$ and denote their elements simply by $(a_1,a_2,\ldots,a_n)$ whenever it is clear from the context whether we are in $SO(2n+1)$ or $SO(2n)$ case. We are using the convention that the exponential map $exp:\lie{t}_{SO(2)} \rightarrow T_{SO(2)}$ is
given by $L(a)\rightarrow R(2\pi a),$ that is $S^1\cong \bb{R}/\bb{Z}$. 
Moreover we choose the positive Weyl chambers to consist of matrices with 
 $a_1\geq  a_2\geq a_3\geq \ldots\geq a_n\geq 0$ in the case $G=SO(2n+1)$, and $a_1\geq a_2\geq a_3\geq \ldots\geq a_{n-1}\geq|a_n|$ in the case $G=SO(2n)$. 

Note that the real dimension of the coadjoint orbit through a point $\lambda=(\lambda_1,\ldots, \lambda_n)$ with
$$\lambda_1=\ldots = \lambda_{n_1} >\lambda_{n_1+1}=\ldots =\lambda_{n_2}> \ldots >\lambda_{n_{m}+1}=\ldots =\lambda_n,$$ 
and $\lambda_n \neq 0$ is equal to
$$\dim \cc{O}_{\lambda}=\begin{cases} 2n^2-(\,k_1(k_1-1) + \ldots + k_{m+1}(k_{m+1}-1)\,) & \textrm{ if }G=SO(2n+1),\\
2n(n-1)-(\,k_1(k_1-1) + \ldots + k_{m+1}(k_{m+1}-1) \,)& \textrm{ if }G=SO(2n),\end{cases}$$
where $k_j=n_{j+1}-n_j$, $n_{m+1}=n$. If $ \lambda_n = 0$ one needs to subtract $\frac 1 2  k_{m+1}(k_{m+1}+1)$ in the $SO(2n+1)$ case, or subtract $\frac 1 2  k_{m+1}(k_{m+1}-1)$ in the $SO(2n)$ case.

The root system of the group $SO(2n+1)$ consists of vectors $\pm e_j$, $j=1, \ldots n$, of squared length $1$, and of vectors $\pm (e_j\pm e_k)$, $j\neq k$, of squared length $2$ in the Lie algebra dual $\lie{t}_{SO(2n+1)}^*\cong \bb{R}^n$. 
Therefore this root system for $SO(2n+1)$ is non-simply laced. 
Note that $$\left\langle  (e_j\pm e_k)^{\vee} , \lambda \right\rangle= 2 \frac{\left\langle e_j\pm e_k,\lambda\right\rangle}{\left\langle e_j\pm e_k,e_j\pm e_k\right\rangle}=\lambda_j \pm \lambda_k$$ and 
$$\left\langle (e_j)^{\vee} , \lambda \right\rangle=2 \frac{\left\langle e_j,\lambda\right\rangle}{\left\langle e_j,e_j\right\rangle}=2\lambda_j.$$

The root system for $SO(2n)$ is simply laced and consists of vectors  $\pm (e_j\pm e_k)$, $j\neq k$, of squared length $2$. 
Note that $$\left\langle  (e_j\pm e_k)^{\vee} , \lambda \right\rangle= 2 \frac{\left\langle e_j\pm e_k,\lambda \right\rangle}{\left\langle e_j\pm e_k,e_j\pm e_k\right\rangle}=\lambda_j \pm \lambda_k.$$

Recall that $r=r_G(\lambda)=\min \{ \langle \alpha_j^{\vee},\lambda \rangle; \alpha_j^{\vee} \textrm{ a coroot, }\langle \alpha_j^{\vee},\lambda \rangle>0\}$. Using the above analysis of the root systems we can calculate that for $\lambda$ in the positive Weyl chamber
$$r_{SO(2n+1)}(\lambda)=\begin{cases}\min\{\lambda_{n_1}-\lambda_{n_1+1}, \lambda_{n_2}-\lambda_{n_2+1},\ldots, \lambda_{n_{m}}-\lambda_{n},\, 2 \lambda_n\} & \textrm{ if }\lambda_n\neq 0,\\
\min\{\lambda_{n_1}-\lambda_{n_1+1}, \lambda_{n_2}-\lambda_{n_2+1},\ldots, \lambda_{n_{m}}\} & \textrm{ if }\lambda_n= 0\end{cases}$$
$$ r_{SO(2n)}(\lambda)=\min\{\lambda_{n_1}-\lambda_{n_1+1}, \lambda_{n_2}-\lambda_{n_2+1},\ldots, \lambda_{n_m}-\lambda_{n_{m+1}}, \lambda_{n_m}+\lambda_{n_{m+1}}\}.$$

Now we are ready to state the main theorem about the Gromov width of coadjoint orbits of the special orthogonal group.
\begin{theorem}\label{mainso}
Consider the coadjoint orbit $\cc{O}_{\lambda}$ of the special orthogonal group passing through a point $\lambda=(\lambda_1,\ldots,\lambda_n) \in \,\lie{t}^*_+$ in the positive Weyl chamber (chosen above)
$$\lambda_1=\ldots = \lambda_{n_1} >\lambda_{n_1+1}=\ldots =\lambda_{n_2}> \ldots >\lambda_{n_{m}+1}=\ldots =\lambda_n,$$ 
satisfying a condition
$$ (*)\,\,\,\,(\lambda_n \neq \lambda_{n-1})\, \vee \,(\lambda_n=0)\, \vee \,(\lambda_n \geq r_G(\lambda)).$$
The Gromov width of $\cc{O}_{\lambda}$ is at least $r_G(\lambda)$, that is,
$$r_{SO(2n+1)}(\lambda)=\begin{cases}\min\{\lambda_{n_1}-\lambda_{n_1+1}, \lambda_{n_2}-\lambda_{n_2+1},\ldots, \lambda_{n-1}-\lambda_{n},\, 2 \lambda_n\} & \textrm{ if }\lambda_n\neq 0,\\
\min\{\lambda_{n_1}-\lambda_{n_1+1}, \lambda_{n_2}-\lambda_{n_2+1},\ldots, \lambda_{n-1}\} & \textrm{ if }\lambda_n= 0\end{cases}$$
if $G=SO(2n+1)$, and  
$$ r_{SO(2n)}(\lambda)=\min\{\lambda_{n_1}-\lambda_{n_1+1}, \lambda_{n_2}-\lambda_{n_2+1},\ldots, \lambda_{n_m}-\lambda_{n_{m+1}}, \lambda_{n_m}+\lambda_{n_{m+1}}\}.$$if $G=SO(2n).$ For orbits $\cc{O}_{\lambda}$ with $\lambda$ that do not satisfy the condition $(*)$ their Gromov width is at least $\lambda_n=\min\{r_G(\lambda), \lambda_n\} .$
\end{theorem}
\subsection{Organization of the paper and acknowledgements.} Section \ref{tech} provides background about Hamiltonian actions and technical ingredients while Section \ref{gtaction} briefly reviews the Gelfand-Tsetlin action in the general setting. The main result is proved separately for the unitary group (Section \ref{un}) and for the special orthogonal groups (Section \ref{son}).

The author is very grateful to Yael Karshon for suggesting this problem and helpful conversations during my work on this project. 
The author also would like to thank Tara Holm and Alexander Caviedes Castro for useful discussions.

\section{Technical ingredients.}\label{tech}
We obtain technical ingredients needed here by generalizing a result of Lisa Traynor \cite[Proposition 5.2]{T} to Proposition \ref{fromlt} and applying the ``classification" result of Karshon and Lerman (\cite{KL}).
While editing this manuscript we found out that equivalent to Proposition \ref{fromlt} was already proved in \cite[Lemma 3.11]{Sch}. It was used by Guangcun Lu in \cite{Lu} to prove one of the claims in his Proposition 1.3. This claim is almost equivalent to our Proposition \ref{anypreimage} (note that Lu uses a different normalization convention). We are still presenting our proof here for completeness, and to show the relation with the result of Lisa Traynor. However we point out that the independent work of Guangcun Lu, \cite{Lu}, was published before our work, and we encourage the reader to consult this reference. The reader familiar with Proposition 1.3 in \cite{Lu} may go directly to Section \ref{gtaction}.

Let $2N$ be the dimension of $\cc{O}_{\lambda}$.
We start with generalizing a result by Lisa Traynor \cite[Proposition 5.2]{T}.
Define the following subsets of $\bb{R}^N$
$$ \Box^N (\pi):=\{0<x_1,\ldots,x_N <\pi\},$$
$$\triangle^N (r):=\{0<y_1,\ldots, y_N;\; y_1+\ldots +y_N<r\}.$$
Equip their product $\Box^N (\pi) \times \triangle^N (r)$ with the symplectic form induced from the standard symplectic structure on $\bb{R}^{2N}$, namely $\sum dx_j\wedge dy_j$.
When the dimension $N$ is understood, we simply write $\Box (\pi)$ and $\triangle (r)$. 
Throughout the paper we use 
$$B^{K}_{\pi\,r}=B^{K}(r)=\{x \in \bb{R}^K;\;|x|^2 <r\}$$
to denote an open $K$-dimensional ball of radius $\sqrt{r}$, i.e. of capacity $\pi \,r$. 
(Note that in \cite{T} $B$ denotes closed balls).

\begin{proposition}\label{fromlt}
 For any $\rho <r$ there is a symplectic embedding of $2N$-dimensional ball $B^{2N}({\rho})$ of radius $\sqrt{\rho}$ (i.e. of capacity $\pi \rho$) into $\Box^N (\pi) \times \triangle^N (r)$. (Both sets are considered as subsets of $\bb{R}^{2N}$ with the standard symplectic form $\sum dx_j\wedge dy_j$.)
\end{proposition}
\begin{proof} 
 
There is a symplectic embedding $\Psi \colon \Box^N (\pi) \times \triangle^N (r) \rightarrow \,B^{2N}(r)$ into the open ball of radius $\sqrt{r}$
\begin{equation*}\label{stdembedding}
\Psi(x_1,y_1, \ldots,x_N,y_N)=(\sqrt{y_1} \cos(2x_1),-\sqrt{y_1} \sin(2x_1), \ldots, \sqrt{y_N} \cos(2x_N),-\sqrt{y_N} \sin(2x_N) ).
\end{equation*}

Let $SD(r) \subset B^2(r)$ be the {\bf slit disc} radius $\sqrt{r}$:
$$SD(r):=B^2(r)\setminus \{x \geq 0, y=0\}\subset \bb{R}^2.$$
Denote by $SD^N(r)$ the corresponding slit polidisc, $SD^N(r):=SD(r)\times \ldots \times SD(r)\in \bb{R}^{2N}.$
It is easy to see that
\begin{equation*}
 \Psi(\kt)=\,B^{2N}(r)\cap SD^N(r).
\end{equation*}
Fix any $\rho <r$ and choose any area preserving diffeomorphism (so also preserving symplectic form)
$$\sigma^{\rho}\colon B^2(\rho) \rightarrow \,Im\, \sigma^{\rho} \subset SD(\rho + \frac 1 N (r -\rho) )\subset SD(r)$$
such that if $x^2+y^2 \leq a$ then $|\sigma^{\rho}(x,y)|^2 \leq a + \frac 1 N (r -\rho)$. Let $\Psi^{\rho}$ be the ``product'' of $N$ $\sigma^{\rho}$'s:
$$\Psi^{\rho} \colon  B^2(\rho)\times \ldots \times  B^2(\rho)\, \rightarrow \,SD^N(\rho + \frac 1 N (r -\rho) ),$$
$$\Psi^{\rho} (x_1,y_1, \ldots,x_N,y_N)=(\sigma^{\rho}(x_1,y_1), \ldots, \sigma^{\rho}(x_N,y_N)).$$
The map $\Psi^{\rho}$ is symplectic as a product of symplectic maps. Furthermore,
$$\Psi^{\rho}(B^{2N}(\rho))\subset \,B^{2N}(r)\cap SD^N(r)=\Psi(\kt),$$
because if $\sum (x_i^2+y_i^2) <\rho$ then $\sum |\sigma^{\rho}(x_i,y_i)|^2 < \rho + N \,\frac 1 N (r -\rho)=r$.
Therefore $\Psi \inv \circ \Psi^{\rho}$ gives symplectic embedding of $B^{2N}(\rho)$ into $\kt$.
\end{proof}

\begin{corollary} If there is a symplectic embedding $\kt \hookrightarrow (M,\omega)$, then 
the Gromov width of $M$ is at least $\pi r$, because for any $\rho <r$ we have a symplectic embedding $B_{\pi \rho} \hookrightarrow M$. 
\end{corollary}
Therefore to find lower bounds for Gromov width, instead of looking for embeddings of symplectic balls, we can look for embeddings of $\kt$. This is exactly how we will proceed. First we review some properties of momentum maps.

An effective action of a torus $T$ on a symplectic manifold $(M, \omega_M) $ is called a \textbf{Hamiltonian action} if there exists a $T$-invariant map
$\Phi \colon M \to \lie{t}^*$, called the \textbf{momentum map} (or moment map), such that
\begin{equation} \label{def moment}
        \iota(\xi_M) \omega =  d \left< \Phi,\xi \right>
\quad \forall \ \xi \in \lie{t},
\end{equation}
where $\xi_M$ is the vector field on $M$ generated by $\xi \in \lie{t}$. Then $M$ is referred to as a {\bf Hamiltonian $T$ manifold}. 
The spaces $\lie{t}^*$, $\lie{t}$ and $\bb{R}^{\dim T}$ are isomorphic though not canonically. Once a specific isomorphism is chosen one can view a momentum map as a map to $\bb{R}^{\dim T}$. Throughout the paper we identify $Lie(S^1)$ with $\bb{R}$ using the convention that the exponential map $\exp \colon Lie(S^1) \rightarrow \bb{R}$ is given by $t \rightarrow e^{2 \pi i t}$ ( so $S^1 \cong \bb{R} /\bb{Z}$).

If $\dim T=\frac 1 2 \dim M$ the Hamiltonian action is called toric. We call $M$ a {\bf proper Hamiltonian $T$ space} if there exists an open and convex subset $\cc{T} \subset \lie{t}^*$, containing $\Phi(M)$ and such that the moment map $\Phi: M \rightarrow \cc{T}$ is proper as a map to $\cc{T}$.
In particular if $M$ is compact, then it is also proper.

\begin{example} Consider the $T^2$ action on $\bb{C}^2$ by $$(e^{it_1},e^{it_2})\dot (z_1,z_2)=(e^{it_1}z_1,e^{it_2}z_2).$$
Then $\Phi \colon \bb{C}^2 \colon \rightarrow \bb{R}^2$ given by 
$$\Phi (z_1,z_2)=(-\pi |z_1|^2,-\pi |z_2|^2)$$
is a momentum map. The image of the momentum map is the (closed) third orthant and the image of a ball of capacity $a$ (so of radius $\sqrt{\frac a \pi}$) is presented in the figure below. 
\begin{center}
\includegraphics[width=0.25\textwidth]{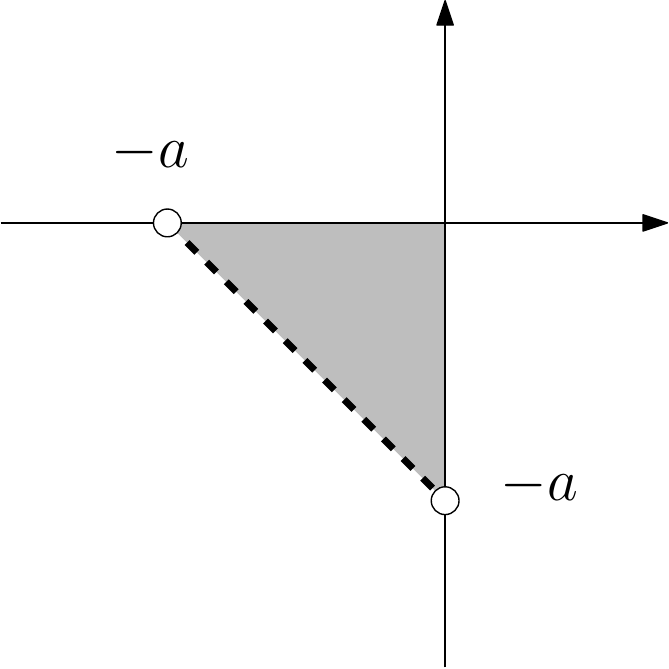}
\end{center}
\end{example}
Compact connected toric manifolds are classified by their momentum map image (Delzant Theorem). Karshon and Lerman generalized this theorem to the case of  non compact manifolds with proper momentum map (\cite{KL}). Using their work we conclude the following proposition.
\begin{proposition}\label{stdpreimage}
 For a connected, proper Hamiltonian $T^N$ space $(M^{2N}, \omega_M)$, with a momentum map $\Phi \colon M \to \lie{t}^*$, and for a subset $S \subset \,int\,\Phi(M)$, $S=W(\triangle^N(r))$ for some $W \in GL(n,\bb{Z})$, we have that $(\Phi \inv (S), \omega_M)$ is symplectomorphic to $T^N \times S$ with the symplectic form $\omega$ given by
$$\omega_{(p,q)}(v,\xi)=\xi(v)=-\omega_{(p,q)}(\xi,v),\;\;\;\omega_{(p,q)}(v,v')=0=\omega_{(p,q)}(\xi,\xi'),$$
for any $(p,q) \in T^N \times S$, $\xi, \xi' \in T_{q}S=\lie{t}^*$, and any $v,v' \in T_pT^N=\lie{t}$.
\end{proposition}
\begin{proof}
The space  $T^N \times S$ with the above symplectic form and the $T^N$ action is on the first factor via the group multiplication is a connected, proper Hamiltonian $T^N$ manifold with momentum map image $S$ (contractible). The space $(\Phi \inv (S), \omega_M)$ is also a connected, proper Hamiltonian $T^N$ manifold. According to Proposition 3.9 in \cite{KL} such manifolds are classified by their momentum map image. Therefore $T^N \times S$ and $\Phi \inv (S)$ are symplectomorphic. Note that the above theorem is true also for more general $S$, however for our purposes it is enough to consider only $S=W(\triangle^N(r))$.
\end{proof}
Recall that we work with the convention $ S^1 =\bb{R}/\bb{Z}.$
Therefore there exists a symplectic embedding
$$(\,(0,1)^N \times S, \sum dx_j \wedge dy_j) \hookrightarrow (T^N \times S, \sum dx_j \wedge dy_j)\cong (\Phi \inv (S), \omega_M).$$
\begin{proposition}\label{anypreimage}
 If $W(\triangle(r)) \subset \,int\,\Phi(M)$, for some $W$ with $\pm W \in  SL(N;\bb{Z})$, then for any $\rho<r$ a ball $B^{2N}_{\rho}=B^{2N}(\rho/\pi)$ of capacity $\rho$ embeds symplectically into $(M,\omega_M)$, and thus the Gromov width of $M$ is at least $r$. 
\end{proposition}
\begin{proof}
First suppose that $W \in  SL(N;\bb{Z})$. According to Proposition \ref{stdpreimage}, 
$$\Phi \inv (\,W(\triangle^N(r))\,)\cong (T^N\times W(\triangle^N(r))\,,\omega) \supset ((0,1)^N \times W(\triangle^N(r)), \sum dx_j \wedge dy_j).$$ 
 Notice that the map $$(Id,W) \colon ((0,1)^N \times \triangle^N(r), \sum dx_j \wedge dy_j)\,\rightarrow ((0,1)^N \times W(\triangle^N(r)),  \sum dx_j \wedge dy_j) $$
is a symplectomophism because $\det (Id,W)=\det W=1.$ 
Also 
$$\left((0,1)^N \times \triangle(r), \,\sum dx_j \wedge dy_j\right)\cong \left((0,\pi)^N \times \triangle(r/\pi),\,\sum dx_j \wedge dy_j\right)$$
are symplectomorphic via $(x,y) \rightarrow (\pi x, y/\pi)$.
Therefore $ \left((0,\pi)^N \times \triangle(r/\pi),\,\sum dx_j \wedge dy_j\right)$ can be symplectically embedded into $\left(\Phi \inv (\,W(\triangle(r)))\, , \; \omega_M\right)$. Together with Proposition \ref{fromlt} this gives that for any $ \rho <r$, a ball of capacity $\rho$, $(B^{2N}(\rho/\pi),  \sum dx_j \wedge dy_j)$, can be symplectically embedded into $M$. If $-W \in  SL(N;\bb{Z})$ then we obtain a symplectic embedding of $(B^{2N}(\rho/\pi), - \sum dx_j \wedge dy_j)$ into $M$, but $(B^{2N}(\rho/\pi), - \sum dx_j \wedge dy_j)$ and $(B^{2N}(\rho/\pi),  \sum dx_j \wedge dy_j)$ are symplectomorphic.
\end{proof}

\section{Gelfand-Tsetlin torus action.}\label{gtaction}
In this Subsection we describe the Gelfand-Tsetlin (sometimes spelled Gelfand-Cetlin, or Gelfand-Zetlin) system of action coordinates, which originally appeared in \cite{GS1}. It is related to the classical Gelfand-Tsetlin polytope introduced in \cite{GTs}. Here we only briefly recall necessary facts about this action and refer the reader to \cite{GS1,NNU,Pun,Pso,Pthesis,K}

Let $G$ be a compact, connected Lie group and $\cc{O}_{\lambda}$ its coadjoint orbit. Consider a sequence of subgroups $G=G_k \supset G_{k-1}\supset \ldots \supset G_1.$  Inclusion of $G_j$ into $G$ gives an action of $G_j$ on $\cc{O}_{\lambda}$. This action is Hamiltonian with momentum map $\Phi^j$, where $\Phi^j$ is the composition of the $G$-momentum map $\Phi$ and the projection $p_j:\lie{g}^* \rightarrow \lie{g}_j^*$. Choose maximal tori, $T_{G_j}$, and positive Weyl chambers for each group $G_j$ in the sequence. Every $G_j$ orbit intersects the positive Weyl chamber $(\lie{t}_{G_j})^*_+$ exactly once. This defines a continuous (but not everywhere smooth)  map $s_j:\lie{g}_j^*\rightarrow (\lie{t}_{G_j})^*_+$. Let $\Lambda^{(j)}=(\lambda^{(j)}_1,\ldots,\lambda^{(j)}_{rk\,G_j})$ denote the composition $s_j \circ \Phi^j$:
\begin{displaymath}
\xymatrix{
\cc{O}_{\lambda} \ar[r]^{\Phi^j} \ar[rd]_{\Lambda^{(j)}} &
\lie{g}_j^*  \ar[d]^{s_j}&
 \\
&(\lie{t}_{G_j})^*_+
} 
\end{displaymath}
The functions $\{\Lambda^{(j)}\}$, $j=1, \ldots, k-1$, form the
\textbf{ Gelfand-Tsetlin system} which we denote by $\Lambda: \cc{O}_{\lambda} \rightarrow \bb{R}^{n(n-1)/2}$. Let $U$ denote the subset of $\cc{O}_{\lambda}$ on which the Gelfand-Tsetlin functions do not coincide ``unnecessarily"
$$U =\{A \in \cc{O}_{\lambda};\,\lambda^{(j)}_k(A)=\lambda^{(j)}_{k+1}(A) \textrm{ if and only if }\lambda^{(j)}_k=\lambda^{(j)}_{k+1} \textrm{ on the whole }\cc{O}_{\lambda}\}.$$
The Gelfand-Tsetlin functions have many useful properties. The ones we are interested in are summarized in the following proposition (for more details see for example \cite{GS1, Pthesis}).
\begin{proposition}\label{gtsummary} In the case of a coadjoint action of $G=U(n)$, $SO(2n+1)$ or $SO(2n)$ on an orbit $\cc{O}_{\lambda}$ through $\lambda \in (\lie{t}_G)^*_+$ (and appropriately chosen sequences of subgroups) the Gelfand-Tsetlin functions are smooth on the open dense subset $U \subset \cc{O}_{\lambda} $ defined above. Moreover, $U$ is equipped with a Hamiltonian action of a torus $T_{GT}$, called the {\bf Gelfand-Tsetlin torus}, of dimension equal to the complex dimension of $\cc{O}_{\lambda} $. This action makes $U$ into a proper toric manifold. The momentum map consists of those coordinates of $\Lambda_{|U}$ which are not constant on the whole orbit. The closure of the momentum map image, $\overline{ \Lambda(U)}$, is the Gelfand-Tsetlin polytope $\cc{P}$ (defined carefully below). In particular $\Lambda^{-1}(\,int\, \cc{P}) \subset U$.
\end{proposition}

\section{Lower bounds for Gromov width of $U(n)$ coadjoint orbits.}\label{un}

In the case of $G=U(n)$ we apply the above procedure to the sequence of subgroups 
$$U(n) \supset U(n-1) \supset \ldots \supset U(2) \supset U(1).$$
Then the Gelfand-Tsetlin functions at $A \in \cc{O}_{\lambda}$ 
$$\lambda^{(j)}_1(A) \geq  \lambda^{(j)}_2(A) \geq \ldots \geq \lambda^{(j)}_{j}(A),\,\,\,\,j=1, \ldots n-1$$
are the eigenvalues of $j \times j$ top-left submatrix of $A$ ordered in a non-increasing way (due to our choice of positive Weyl chamber).

The classical min-max principle (see for example Chapter I.4 in \cite{CH}) implies that 
\begin{equation*}
 \lambda^{(l+1)}_j(A) \geq  \lambda^{(l)}_j(A) \geq \lambda^{(l+1)}_{j+1}(A).
\end{equation*}
These inequalities cut out a polytope in $\bb{R}^{n(n-1)/2}$, which we denoted by $\cc{P}$, and  $\Lambda(\cc{O}_{\lambda})$ is contained in this polytope. In fact, 
$\Lambda(\cc{O}_{\lambda})$ is exactly $\cc{P}$ (\cite{GS1,Pthesis,NNU}). According to Proposition \ref{gtsummary} the number of ``non-trivial "Gelfand-Tsetlin functions, i.e. ones that are not constant on the whole orbit $\cc{O}_{\lambda}$ is equal to $N=\dim_{\bb{C}}\cc{O}_{\lambda}=\dim T_{GT}$. These $N$ functions are the coordinates of momentum map for the $T_{GT}$ action. Therefore $\Lambda(\cc{O}_{\lambda})$ in fact sits in some affine $\bb{R}^N \subset \bb{R}^{n(n-1)/2}$ and one can view the polytope $\cc{P}$ as a polytope in $\lie{t}_{GT}^* \cong \bb{R}^N \subset \bb{R}^{n(n-1)/2}$.

Convenient way to visualize the Gelfand-Tsetlin functions for $U(n)$ is via the standard ladder diagram(\cite{BCKS, NNU}). For chosen $\lambda$,
$$\lambda_1=\ldots = \lambda_{n_1} >\lambda_{n_1+1}=\ldots =\lambda_{n_2}> \ldots >\lambda_{n_{m}+1}=\ldots =\lambda_n,$$
let $Q=Q_{\lambda}$ be an $n \times n$ square with squares $Q_l$ of size $(n_l-n_{l-1})\times (n_l-n_{l-1})$, $l=1,\ldots,m+1$, $n_0=0$, $n_{m+1}=n$, on the diagonal. The {\bf ladder diagram} is the set of boxes below the diagonal squares. Note that the number of boxes in the ladder diagram for $\lambda$ is equal to $N$, the number of non-trivial Gelfand-Tsetlin functions. To refer to particular boxes we think of $Q$ as sitting in the first quadrant of $\bb{R}^2$ and use Cartesian coordinates. Min-max inequalities imply that for every box in the ladder diagram $Q$ its value needs to be between the values of its right neighbor and top neighbor. The coordinates of $T_{GT}$ corresponding to these non-constant functions give an effective Hamiltonian torus action on $U$.

\begin{figure}[h]
		\includegraphics[width=0.4\textwidth]{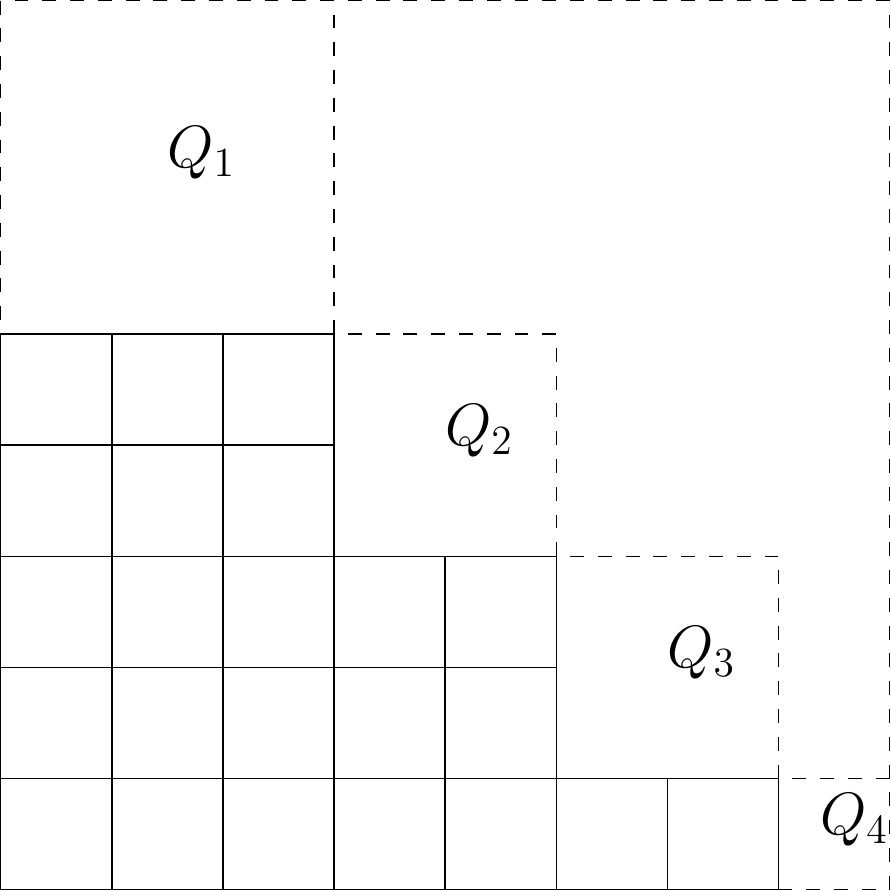}\hspace{16mm}
\includegraphics[width=0.4\textwidth]{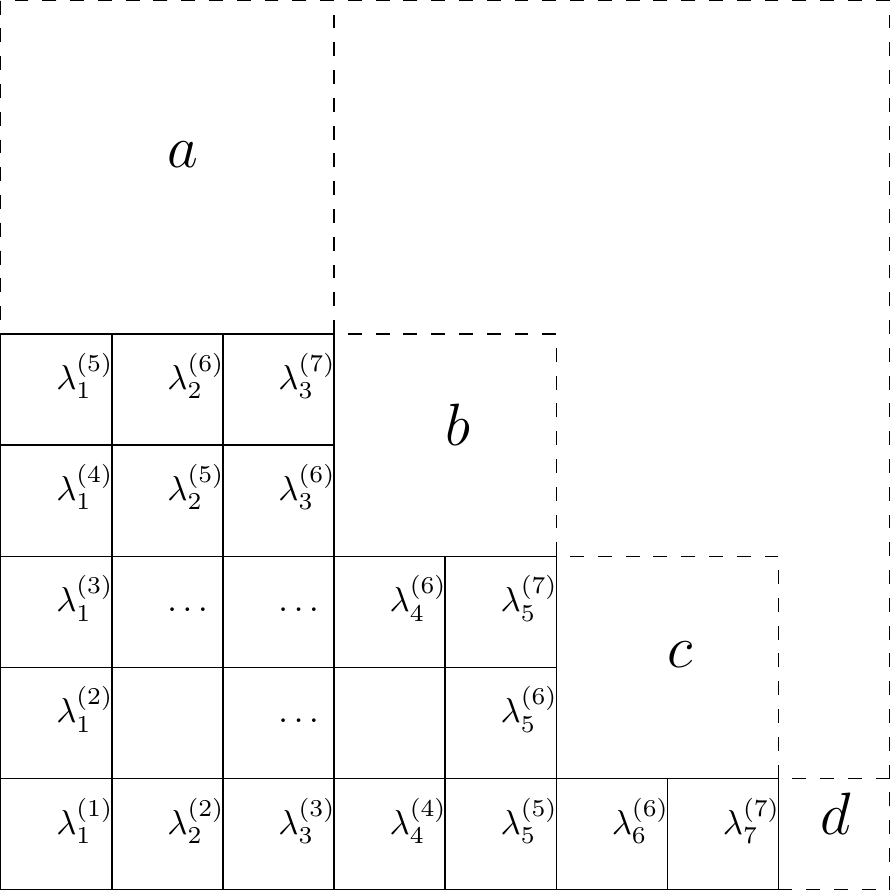}
	\caption{The ladder diagram for $\lambda=(a,a,a,b,b,c,c,d)$ and its filling with the Gelfand-Tsetlin functions.}\label{ladderdiagram}
\end{figure}

For what follows it will be more convenient to index the Gelfand-Tsetlin functions with Cartesian coordinates of the corresponding ladder diagram. Therefore let
$$\lambda_{j,k}:=\lambda^{(j+k-1)}_{j}.$$

We need more careful analysis of the Gelfand-Tsetlin polytope $\cc{P}:=\Lambda(\cc{O}_{\lambda})$. 
Let $\{e_{j,k} \,|\, (j,k) \textrm{ a box in the ladder diagram }\}$ denote the set of generators of $\bb{R}^N\cong \lie{t}^*_{GT}$. As usually, for $x \in \bb{R}^N$ we denote by $x_{j,k}$ its coordinate in $e_{j,k}$ direction. We fix the following ordering of generators to obtain an ordered basis for $\bb{R}^N$: 
$$e_{j,k}\textrm{ proceeds }e_{j',k'}\textrm{ iff }k'>k\textrm{ or }k=k'\textrm{ and }j>j'.$$ \begin{figure}[h]
\includegraphics[width=0.35\textwidth]{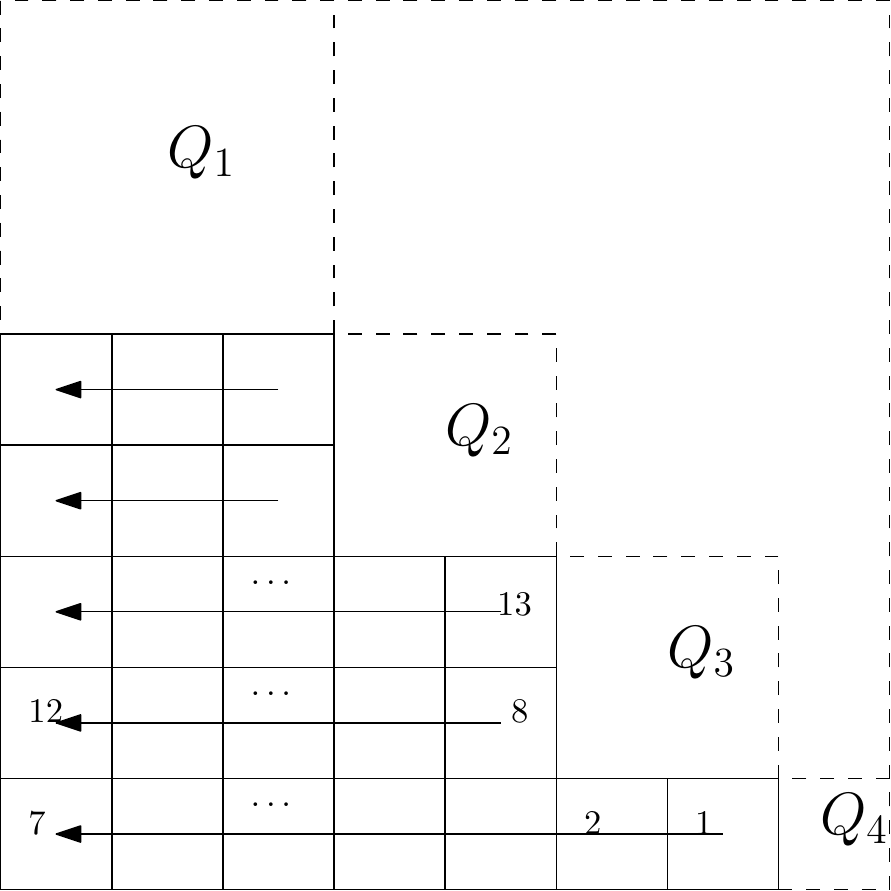}
 	\caption{The ordered basis of $\bb{R}^N$.}\label{orderingbasis}
 \end{figure}

For each $l=1,\ldots, n$ denote by $g(l)$ an integer such that the row $l$ intersects the diagonal box $Q_{n_{g(l)}}$. This implies $$\lambda_{n-l+1}=\lambda_{n_{g(l)}}.$$
For example, in the situation presented on Figure \ref{ladderdiagram} we have $g(5)=g(4)=2$, $g(3)=g(2)=3$, $g(1)=4$.

Let  $V=\Lambda( \, \diag(\lambda_n,\ldots,\lambda_1))\in \cc{P}$ be the image of $\diag(\lambda_n,\ldots,\lambda_1) \in \lie{u}(n)^*$,  that is $V =\sum\,V_{j,k}\,e_{j,k}$ where
$$V_{j,k}=\lambda_{j,k}( \, \diag(\lambda_n,\ldots,\lambda_1)=\lambda_{n-k+1}=\lambda_{n_{g(k)}}.$$ Note that $V$ is a vertex of $\cc{P}$ because all Gelfand-Tsetlin functions when evaluated at $ \diag(\lambda_n,\ldots,\lambda_1)$ are equal to their lower bounds. The vertex $V$ does not need to be smooth: there might be more than $N$ edges in $\cc{P}$ starting from $V$.

\begin{figure}[h]
\includegraphics[width=0.4\textwidth]{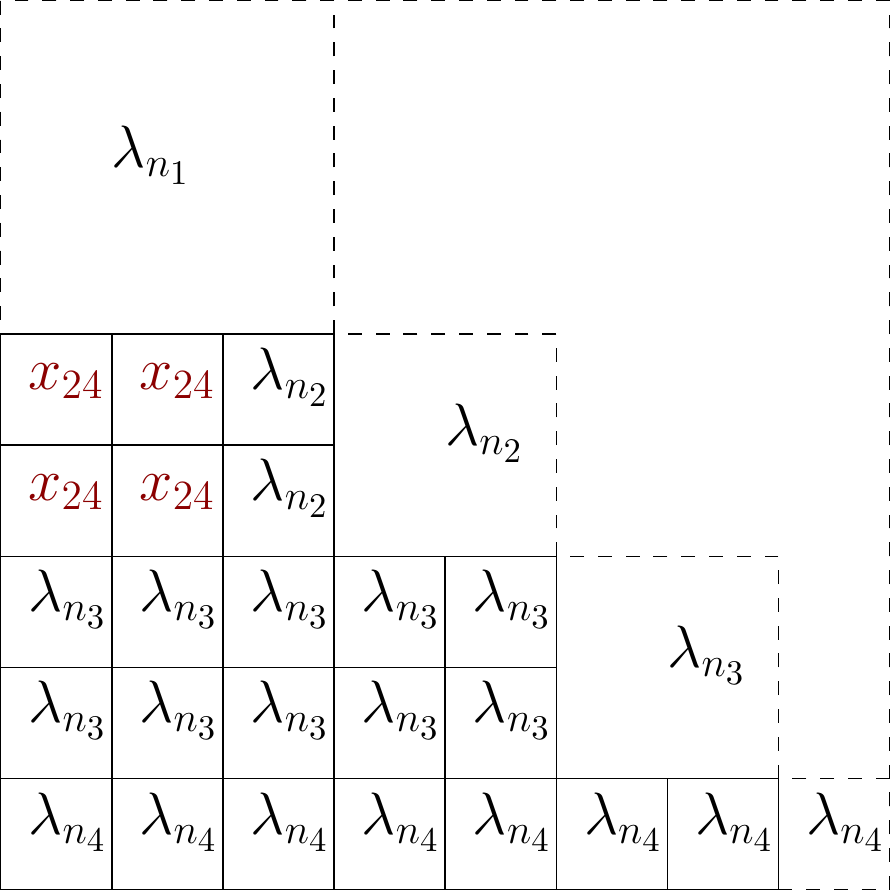}\hspace{16mm}
\includegraphics[width=0.4\textwidth]{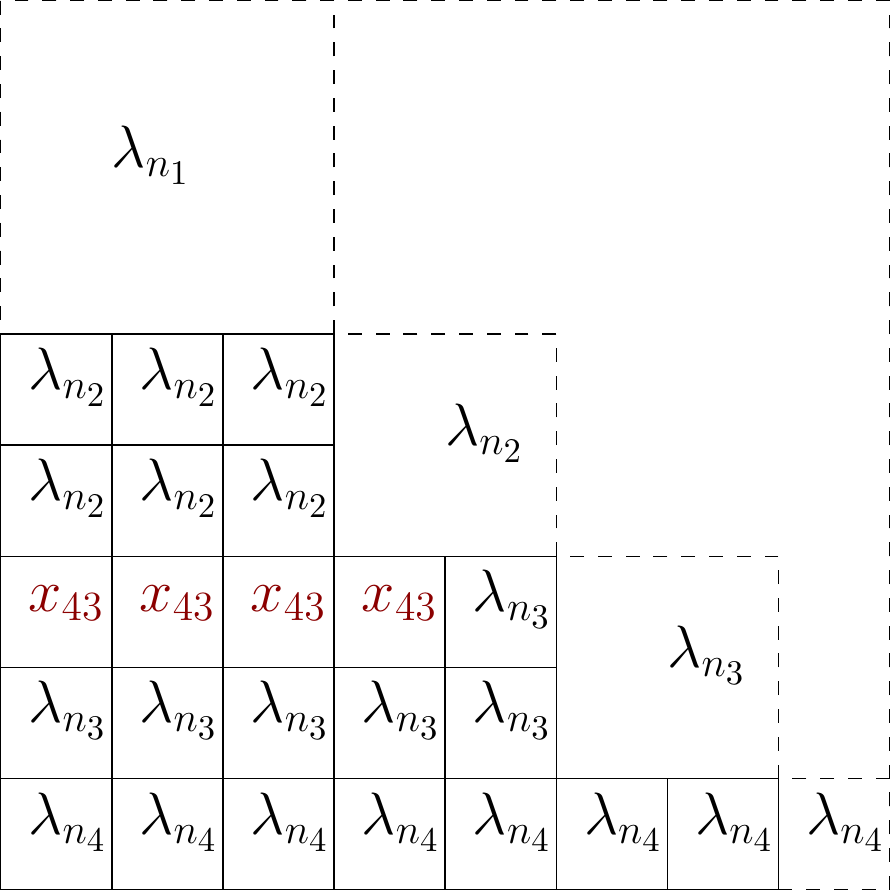}
	\caption{Elements of sets $E_{2,4}$, ($x_{24} \in [\lambda_{n_2},\lambda_{n_1}]$) and $E_{4,3}$, ($x_{43} \in [\lambda_{n_3},\lambda_{n_2}]$).}\label{edges}
\end{figure}

Take $(s,l)$ such that $(s,l)$-th box is in the ladder diagram. Define the following $N$ subset of $\cc{P}$,
\begin{align*}
E_{s,l}=\{x \in \bb{R}^N\,|\,x_{s,l} &\in [\lambda_{n_{g(l)}}, \lambda_{n_{g(l)-1}}],\;\; \\
x_{j,k}&=x_{s,l} \textrm{ for }j\leq s,\; l \leq k  \leq n-n_{g(l)-1},\;\;\\
x_{j,k}&=\lambda_{n-k+1} \textrm{ for other }(j,k)\;\}.
\end{align*}
In other words, $E_{s,l}$ is the set of points with almost all coordinates equal to the coordinates of $V$. We allow the coordinate $x_{s,l}$ to be greater than $V_{s,l}$, and this forces some other coordinates to change as well due to min-max principle. See Figure \ref{edges}.

\begin{lemma}\label{unlengths}
Each subset $E_{s,l}$ is an edge of $\cc{P}$ starting from $V$. The lattice length of $E_{s,l}$ is $\lambda_{n_{g(l)-1}}-\lambda_{n_{g(l)}}$.
\end{lemma}

\begin{proof}
Let $H_{s,l}$ be an affine subspace of $\bb{R}^N$ defined by 
\begin{align*}
x_{j,k}&=x_{s,l} \textrm{ for }j\leq s,\; l \leq k  \leq n-n_{g(l)-1},\\
x_{j,k}&=\lambda_{n_{g(k)}} \textrm{ for other }(j,k)\;.
\end{align*}
Then $E_{s,l}=H_{s,l} \cap \cc{P}$ is a face of $\cc{P}$ because $\cc{P}$ is fully contained in the closure of one component of $\bb{R}^N \setminus H_{s,l}$. 
Any $x \in H_{s,l} \cap \cc{P}$ is determined by the value of $x_{s,l} \in (\lambda_{n_{g(l)}}, \lambda_{n_{g(l)-1}}),$ therefore this face is $1$-dimensional and of lattice length $\lambda_{n_{g(l)-1}}-\lambda_{n_{g(l)}}$.  Each of the $N$ boxes $(s,l)$ in the ladder diagrams gives such an edge.
Moreover $V$ belongs to each $E_{s,l}$.
\end{proof}

Let $w_{s,l}$ denote the primitive vector in the direction of $E_{s,l}$ (starting from $V$), that is
$$w_{s,l}=\sum e_{j,k},$$
where the sum is over $j,k$ such that $j\leq s,\; l \leq k  \leq n-n_{g(l)-1}$.
Recall that
$$r=r_{U(n)}(\lambda)=\min\{\,\lambda_{n_1}-\lambda_{n_1+1}, \lambda_{n_2}-\lambda_{n_2+1},\ldots, \lambda_{n_m}-\lambda_{n_m+1} \}.$$ 
Convexity of $\cc{P}$ and Lemma \ref{unlengths} imply that the following set is a subset of $\cc{P}$
$$R:=\textrm{convex hull }\{V,V+r\,w_{s,l}\,|\, (s,l)\textrm{ a box in the ladder diagram }\}.$$

\begin{lemma}\label{gooduncorner}
 $R$ is $ SL(N,\bb{Z})$-equivalent to the closure of an $N$-dimensional tetrahedron $\triangle^N (r )$,
$$\triangle^N (r )=\{0<y_1,\ldots, y_N;\; y_1+\ldots +y_N<r\}.$$
\end{lemma}
\begin{proof}
Edges of $R$ starting from $V$ are given by $\{r\,w_{s,l}\}.$ Notice that the first non-zero coordinate of $w_{s,l}$ is equal to $1$ and appears on the $e_{s,l}$-th coordinate. Therefore the matrix of vectors
$w_{s,l}$, ordered the same way we ordered the basis elements, is an integral (all entries are $0$ or $1$), lower triangular matrix, with $1$'s on diagonal. Thus it belongs to $SL(N,\bb{Z})$.
\end{proof}
\begin{proof}{\it (of Theorem \ref{width}.) }
 Recall that $\Lambda \inv (Int\,\cc{P}) \subset U$, and  $U$ is a proper Hamiltonian $T_{GT}$-space, $\dim T_{GT}=N= \dim_{\bb{C}}U$. The non-constant coordinates of  $\Lambda|_{U}$ form a momentum map. Lemma \ref{gooduncorner} gives that $W(\triangle^N (r ))=R \subset \,Int\,\cc{P}$ for some change of basis matrix $W \in  SL(N;\bb{Z})$. Therefore, by Proposition \ref{anypreimage} the Gromov width of $\cc{O}_{\lambda}$ is at least $r$, as claimed.
\end{proof}
\section{Lower bounds for Gromov width of $SO(2n)$ and $SO(2n+1)$ coadjoint orbits.}\label{son}

In the case of special orthogonal group $SO(m)$, $m=2n$ or $2n+1$ we work with the following sequence of subgroups
\begin{equation*}
 G_{m}=SO(m) \supset G_{m-1}=SO(m-1) \supset G_{m-2}=SO(m-2) \supset \ldots \supset G_2=SO(2).\end{equation*}
and obtain the following Gelfand-Tsetlin functions. For each $j<m$, 
$$\lambda^{(j)}_1(A) \geq  \lambda^{(j)}_2(A) \geq \ldots \geq \lambda^{(j)}_{\lfloor \frac m 2 \rfloor}(A),\,\,(\geq 0\textrm{ if }j \textrm{ odd })$$ 
are numbers such that $j \times j$ submatrix of $A$ is $SO(j)$ equivalent to $ (\lambda^{(j)}_1(A), \ldots,\lambda^{(j)}_{\lfloor \frac m 2 \rfloor}(A)) \in (\lie{t}_{SO(j)})^*_+.$
Note that for any $\lambda$ the number of the Gelfand-Tsetlin functions that are not constant on the whole orbit $\cc{O}_{\lambda}$ is exactly the complex dimension of the orbit.

Similarly to the unitary case the following inequalities between the above functions need to be satisfied:
$$\begin{cases}
\,\,\,\lambda^{(2k)}_{1} \geq \lambda^{(2k-1)}_{1} \geq \lambda^{(2k)}_{2} \geq \lambda^{(2k-1)}_{2} \geq \ldots \geq \lambda^{(2k)}_{k-1} \geq \lambda^{(2k-1)}_{k-1} \geq |\lambda^{(2k)}_{k}| ,\\
\,\,\,\lambda^{(2k+1)}_{1} \geq \lambda^{(2k)}_{1} \geq \lambda^{(2k+1)}_{2} \geq \lambda^{(2k)}_{2} \geq \ldots \geq \lambda^{(2k+1)}_{k} \geq |\lambda^{(2k)}_{k}| ,\\
\end{cases},$$ for all indices for which it makes sense. Moreover the image $\Lambda(\cc{O}_{\lambda})$ is equal to the polytope $\cc{P}=\cc{P}_{\lambda}$, defined by the above set of inequalities (see \cite{Pthesis}).

\begin{figure}[h]
		\includegraphics[width=1\textwidth]{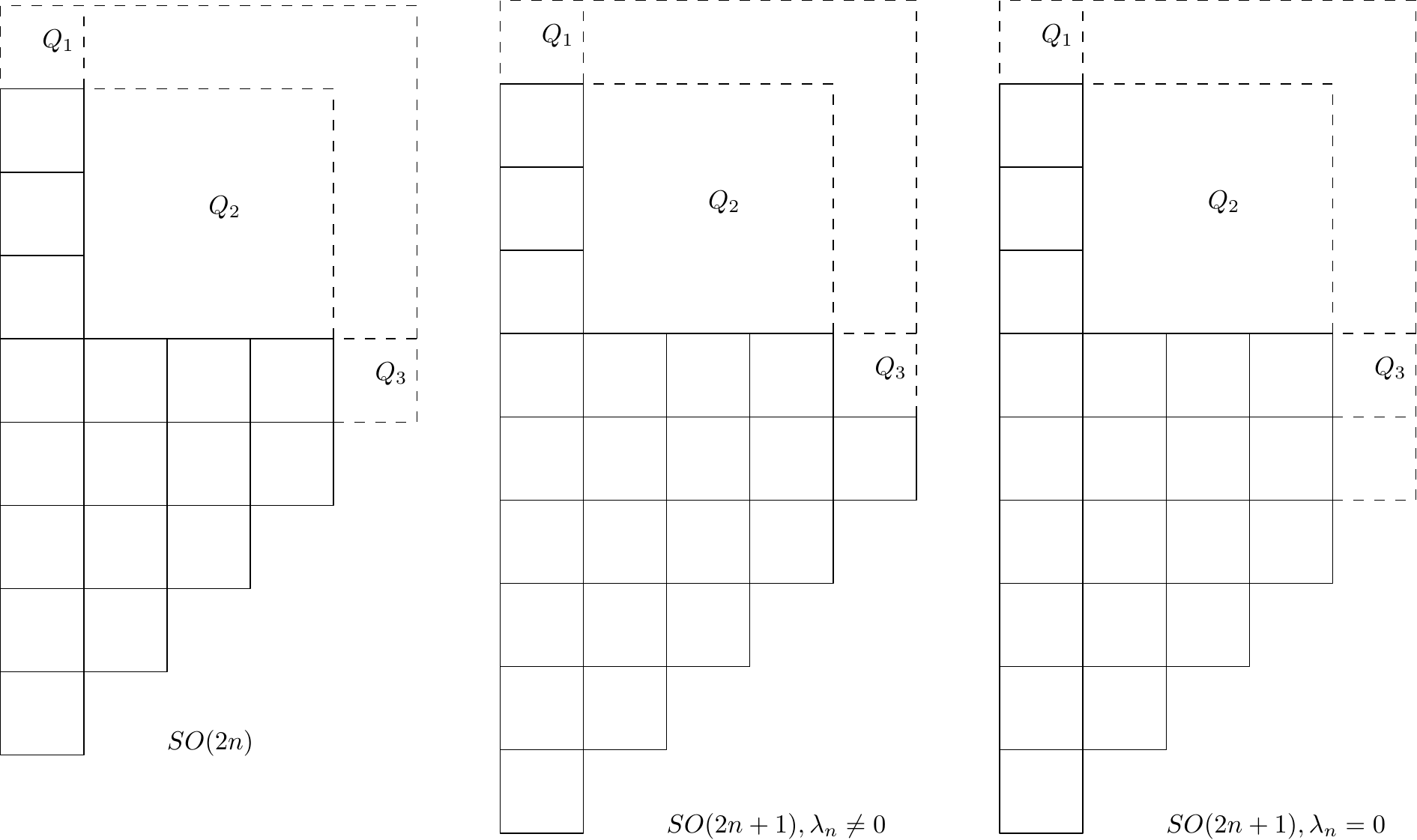}
	\caption{The so-diagrams for $SO(2n+1)$ and $SO(2n)$.}\label{sodiagram}
\end{figure}

To visualize the Gelfand-Tsetlin functions for $G=SO(2n+1)$ and $SO(2n)$ orbits $\cc{O}_{\lambda}$ through $\lambda \in (\lie{t}_G)^*_+$ we add extra boxes from the IV-th quadrant to the $U(n)$-ladder diagram for the corresponding $\lambda$. In the $SO(2n+1)$ case we add $\frac 1 2 (n^2-n)+n$ boxes and in the $SO(2n)$ case add $\frac 1 2 (n^2-n)$ boxes the way presented on the Figure \ref{sodiagram} (first two pictures).

Each diagonal (that is intersection of the diagram with a line of slope $-1$) corresponds to functions with the same superscript. Boxes in diagonal squares $Q_l$ contain functions which are constant on the whole orbit. If $\lambda_n=0$ then all the Gelfand-Tsetlin functions corresponding to the boxes below the last diagonal square are equal to $0$ on the whole orbit. We delete these boxes from the diagram (an example is presented on the third picture in Figure \ref{sodiagram}). Now there is a one-to-one correspondence between the boxes in the diagram and the Gelfand-Tsetlin functions that are not constant on the whole orbit.
Call such a diagram the {\bf so-diagram}. 

The subspace of $(\lie{t}_{G})^*$ spanned by the images of non-constant Gelfand-Tsetlin functions has the dimension equal to the complex dimension of the orbit $\cc{O}_{\lambda}$, which we continue to denote by $N$. We identify this subspace with $\bb{R}^N$, where the basis of $\bb{R}^N$ consists of elements $\{e_{j,k};\,\,(j,
k)\textrm{ is in the so-diagram}\}$ ordered in the following way:
$$e_{j,k}\textrm{ proceeds }e_{j',k'}\textrm{ iff }j<j'\textrm{ or }j=j'\textrm{ and }k>k'.$$\begin{figure}[h]
		\includegraphics[width=0.3\textwidth]{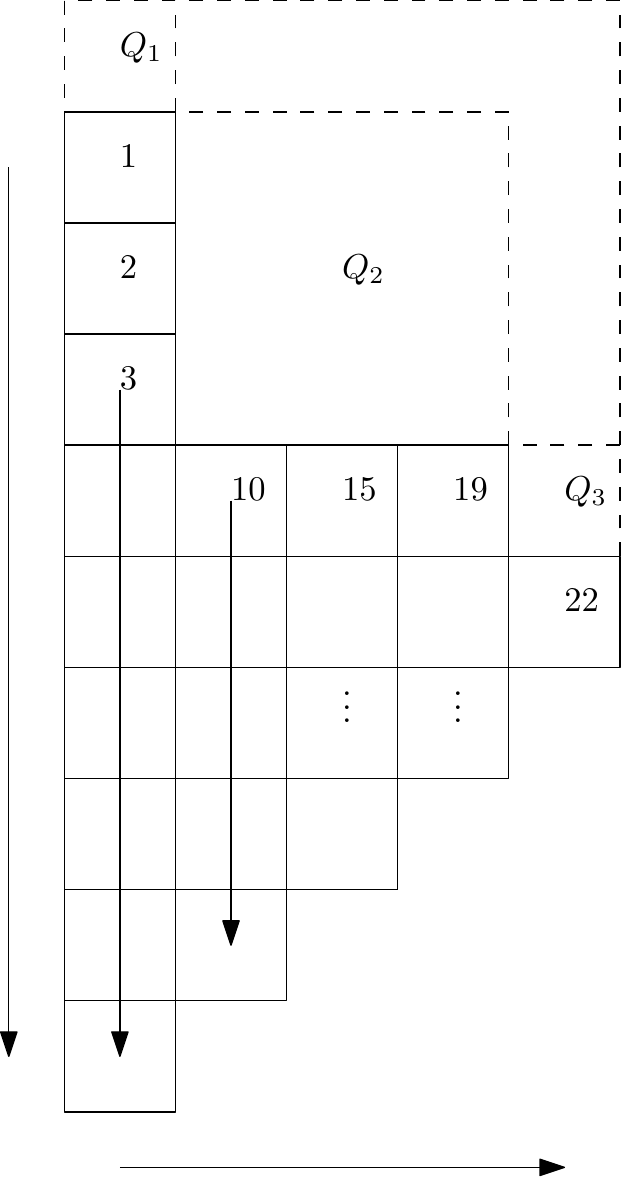}\hspace{16mm}
	\caption{Ordering of the basis.}\label{sodiagrambasis}
\end{figure}

Let $V=\Lambda(\lambda)$. Notice that $V=\sum\,\lambda_j\,e_{j,k}$, with the sum taken over all boxes in the so-diagram. It is the vertex of $\cc{P}$ as all Gelfand-Tsetlin functions attain their maximum at $\lambda$. Recall the definition of $g(j)$ from the previous Section and notice that the column $j$ in the so-diagram intersects the diagonal square $Q_{g(j)}$.

Take $(s,l)$ such that $(s,l)$-th box is in the diagram. 
Define the following affine subspaces of $\bb{R}^N$, one per each box in the so-diagram
\begin{align*}
H_{s,l}=\{x \in \bb{R}^N\,|\,x_{s,l} &\in \bb{R},\;\; \\
x_{j,k}&=x_{s,l} \textrm{ for }k\leq l,\; s \leq j  \leq n_{g(s)},\;\;\\
x_{j,k}&=\lambda_{j} \textrm{ for other }(j,k)\;\}.
 \end{align*} 

Recall that $r=r_G(\lambda)=\min \{ \langle \alpha_j^{\vee},\lambda \rangle; \alpha_j^{\vee} \textrm{ a coroot, }\langle \alpha_j^{\vee},\lambda \rangle>0\}$, that is
$$r_{SO(2n+1)}(\lambda)=\begin{cases}\min\{\lambda_{n_1}-\lambda_{n_1+1}, \lambda_{n_2}-\lambda_{n_2+1},\ldots, \lambda_{n-1}-\lambda_{n},\, 2 \lambda_n\} & \textrm{ if }\lambda_n\neq 0,\\
\min\{\lambda_{n_1}-\lambda_{n_1+1}, \lambda_{n_2}-\lambda_{n_2+1},\ldots, \lambda_{n-1}\} & \textrm{ if }\lambda_n= 0\end{cases}$$
$$ r_{SO(2n)}(\lambda)=\min\{\lambda_{n_1}-\lambda_{n_1+1}, \lambda_{n_2}-\lambda_{n_2+1},\ldots, \lambda_{n_m}-\lambda_{n_{m+1}}, \lambda_{n_m}+\lambda_{n_{m+1}}\}.$$
\begin{lemma}\label{lengths}
For each of these affine hyperplanes the intersection $H_{s,l} \cap \cc{P}$ is an edge of $\cc{P}$. Edges $E_{s,l}$ with $g(s)=m+1$ (that  is $n_{g(s)}=n$) and $s-n < l \leq 0$ in the $SO(2n+1)$ case, $s-n+1 < l \leq 0$ in the $SO(2n)$ case, if they exist, are of lattice length $\lambda_n$. All other edges are of lattice length at least $r_G(\lambda)$.
\end{lemma}
\begin{proof}
For each of these affine hyperplanes, the polytope $\cc{P}$ is contained in a closure of one connected component of $\bb{R}^N \setminus H_{s,l}$. Therefore $E_{s,l}:=H_{s,l}\cap \cc{P}$ is a face of $\cc{P}$. A point in $E_{s,l}$ is uniquely determined by the value $x_{s,l}$. Therefore each $E_{s,l}$ is an edge of $\cc{P}$. The length of these edges depends on the position of $(s,l)$ box, and on whether we are in $SO(2n+1)$ or $SO(2n)$ case. Examples are presented in Figures \ref{oddedges} and \ref{evenedges}. One can easily check that the lengths are given by the following formulas. 
For $SO(2n+1)$
$$ |E_{s,l}|=\begin{cases}\lambda_{n_{g(s)}} & \textrm{ if }s < n_{g(s)} \textrm{ and } s-n < l\leq n_{g(s)}-n ,\\
2\lambda_{n_{g(s)}} & \textrm{ if on the ``rim" i.e. }l=s-n,\\
\lambda_{n_{g(s)}}-\lambda_{n_{g(l)+1}}& \textrm{ otherwise }.\end{cases}$$
For $SO(2n)$
$$ |E_{s,l}|=\begin{cases}  \lambda_{n_{g(s)}} & \textrm{ if } s < n_{g(s)} \textrm{ and }s-n+1 < l\leq n_{g(s)}-n+1 ,\\
2\lambda_{n_{g(s)}} & \textrm{ if on the ``rim" i.e. }l=s-n+1,\\
\lambda_{n_{g(s)}}-|\lambda_{n_{g(l)+1}}| & \textrm{ otherwise (recall: }\lambda_n \textrm{ can be negative)}.\end{cases}$$ \begin{figure}[h]
\includegraphics[width=0.2\textwidth]{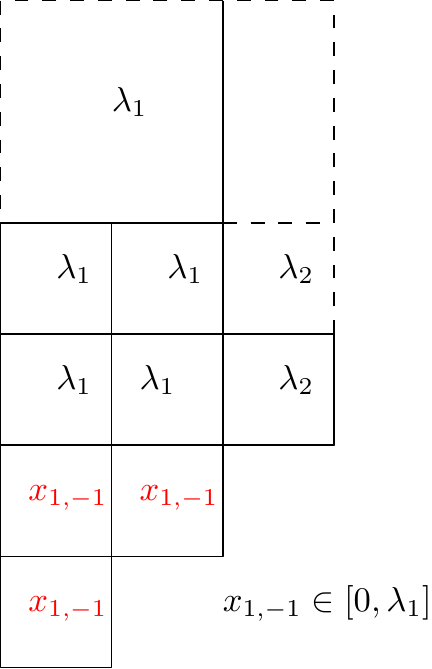}\hspace{10mm}
		\includegraphics[width=0.22\textwidth]{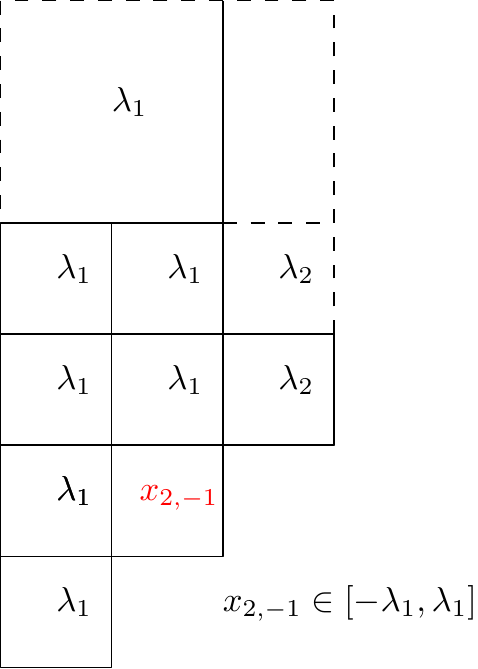}\hspace{10mm}
\includegraphics[width=0.2\textwidth]{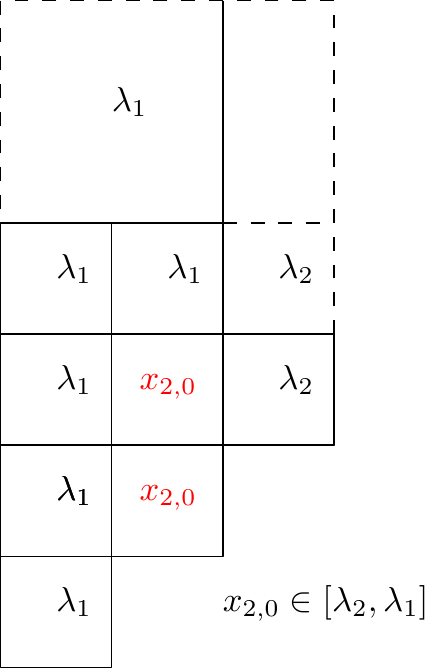}
	\caption{The edges $E_{1,-1}$, $E_{2,-1}$ and $E_{2,0}$ in the polytope for $SO(2n+1)$ coadjoint orbit.}\label{oddedges}
\end{figure}    
\begin{figure}[h]
\includegraphics[width=0.2\textwidth]{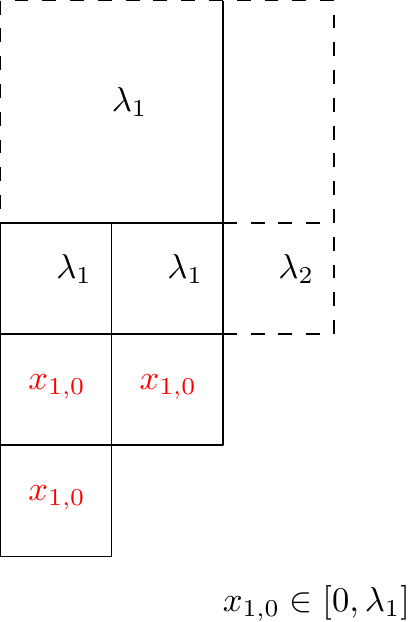}\hspace{10mm}
		\includegraphics[width=0.23\textwidth]{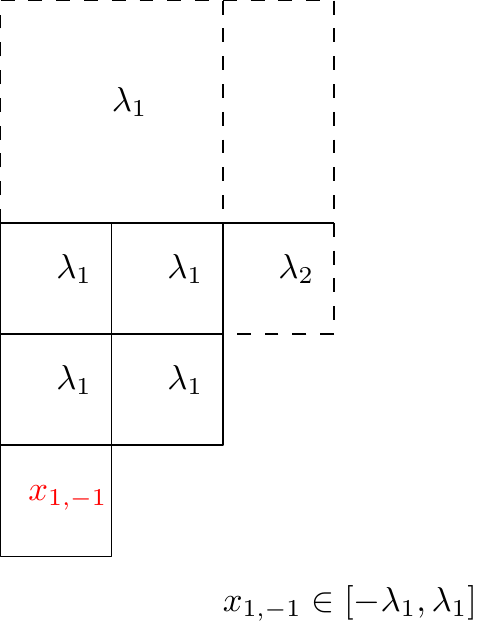}\hspace{10mm}
\includegraphics[width=0.21\textwidth]{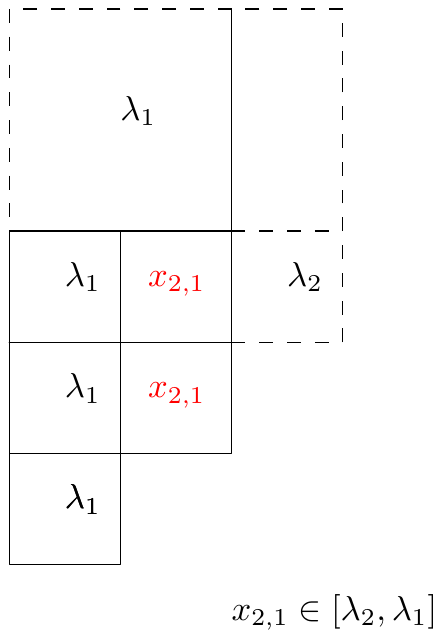}
	\caption{The edges $E_{1,0}$, $E_{1,-1}$ and $E_{2,1}$ in the polytope for $SO(2n)$ coadjoint orbit.}\label{evenedges}
\end{figure} 
All of the above values are positive. The value $\lambda_j$ can be negative if and only if $G=SO(2n)$ and $j=n$. Then $\lambda_n\neq \lambda_{n-1}$. Therefore there are no boxes $(s,l)$ with $g(s)=m+1$.

Note that the values $2\lambda_{n_{g(s)}}$, $\lambda_{n_{g(s)}}-\lambda_{n_{g(l)+1}}$, and $\lambda_{n_{g(s)}}+\lambda_{n_{g(l)+1}}$ are equal to $\left\langle \alpha^{\vee},\lambda \right\rangle$ for some coroot $\alpha^{\vee}$.  
Moreover, for $g(s)\leq m$ 
$$\lambda_{n_{g(s)}} \geq \lambda_{n_{g(s)}}-|\lambda_{n_{g(s)+1}}|=\left\langle \alpha^{\vee},\lambda \right\rangle$$
for some coroot $\alpha^{\vee}$.
Only the edges $E_{s,l}$ with $g(s)=m+1$ (that  is $n_{g(s)}=n$) and $s-n < l \leq 0$ in the $SO(2n+1)$ case, $s-n+1 < l \leq 0$ in the $SO(2n)$ case, if they exist, have lattice length equal to $\lambda_n$ and $\lambda_n$ can be smaller than $r$. 
\end{proof}
\begin{lemma} \label{summarylength} If $\lambda$ satisfies condition $(*)$ then all edges $E_{s,l}$ defined above have lattice length at least $r_G(\lambda)$.
\end{lemma}
\begin{proof}
Recall that the condition $(*)$ says:
$$(\lambda_n \neq \lambda_{n-1})\, \vee \,(\lambda_n=0)\, \vee \,(\lambda_n \geq r_G(\lambda)).$$
If $\lambda_n \geq r_G(\lambda)$ then Lemma \ref{lengths} proves the claim.
If $\lambda_n=0$, or if $\lambda_n \neq \lambda_{n-1}$ and $G=SO(2n)$, then there are no boxes $(s,l)$ with $g(s)=m+1$ and thus, by Lemma \ref{lengths}, all edges have lattice length at least $r_G(\lambda)$.
If  $\lambda_n \neq \lambda_{n-1}$ and $G=SO(2n+1)$ then the only box $(s,l)$ with $g(s)=m+1$ is the $(n,0)$ box, which is on the ``rim" of the so-diagram. The edge $E_{n,0}$ corresponding to this box have lattice lentgh $2\lambda_n$.
\end{proof}
Let $w_{s,l}$ denote the primitive vector in the direction of $E_{s,l}$ (starting from $V$), that is
$$w_{s,l}=\,-\,\sum e_{j,k},$$
where the sum is over $j,k$ such that $k\leq l,\; s \leq j  \leq n_{g(s)}$. Recall that $r=r_G(\lambda)$. Let
$$r':=\begin{cases}r & \textrm{ if $(*)$ is satisfied},\\ \lambda_n &\textrm{ if $(*)$ is not saitisfied (what implies } r>\lambda_n\neq 0)
\end{cases}.$$ 
Denote by
$$R:=\textrm{convex hull }\{V,V+r'\,w_{s,l}\,|\, (s,l)\textrm{ a box in the ladder diagram }\}.$$
Convexity of $\cc{P}$ and Lemmas \ref{lengths}, \ref{summarylength} imply that $R \subset \cc{P}$.
\begin{lemma}\label{goodcorner}
 $R$ is $ \pm SL(N,\bb{Z})$-equivalent to the closure of an $N$-dimensional tetrahedron $\triangle^N (r' )$.
\end{lemma}
\begin{proof}
Edges of $R$ starting from $V$ are given by $\{r\,w_{s,l}\}.$ Notice that the first non-zero coordinate of $w_{s,l}$ is equal to $- 1$ and appears on the $e_{s,l}$-th coordinate. Therefore the matrix of vectors
$w_{s,l}$, ordered the same way we ordered the basis elements, is an integral (all entries are $0$ or $- 1$), lower triangular, with $-1$'s on diagonal. Thus it belongs to $\pm SL(N,\bb{Z})$.
\end{proof}

\begin{proof}{\it (of Theorem \ref{mainso}.) } We need to show that the Gromov width is at least $r'$. The proof is analologous to the unitary case. Subset $U \subset \cc{O}_{\lambda}$ is a proper Hamiltionan $T_{GT}$ space with $\dim T_{GT} = \dim_{\bb{C}} U$. Moreover $\Lambda \inv (R) \subset \Lambda \inv( int\,\cc{P} )\subset U$ and $R = W (\triangle(r' ))$ for some $W \in \pm SL(N,\bb{Z})$. Therefore, by Proposition \ref{anypreimage} the Gromov width of $\cc{O}_{\lambda}$ is at least $r'$, as claimed.
\end{proof}
\begin{remark} Note that if $\lambda$ does not satisfy condition $(*)$, then the Gelfand-Tsetlin polytope $\cc{P}$ is contained between two affine hyperplanes $(x_{s,l}=0)$ and $(x_{s,l}=\lambda_n)$ where $(s,l)=(n-1,0)$ in the $SO(2n+1)$ case and $(s,l)=(n-2,0)$ in the $SO(2n+1)$ case. These hyperplanes are lattice distance $\lambda_n<r$ apart and therefore in that case $\cc{P}$ cannot contain $W(\triangle^N (r ))$ for any $W \in \pm SL(N,\bb{Z})$. This means that $B_r$, a ball of capacity $r$, cannot be equivariantly (with respect to the Gelfand-Tsetlin action) embedded into $U$. There might still exist a symplectic though not equivariant embedding of $B_r$. The Gromov width of this orbit might still be equal to $r$. 
\end{remark}
\bibliographystyle{amsalpha}
\bibliography{GelfandTsetlinRelated.bib}
\end{document}